\documentclass{amsart}
\usepackage{amsmath,amssymb,amsfonts,euscript,graphicx,color}

\newtheorem{pro}[equation]{Proposition}

\newtheorem{lem}[equation]{Lemma}
\newtheorem{ex}[equation]{Example}

\newtheorem{DEF}[equation]{Definition}
\newtheorem{rem}[equation]{Remark}

\def\andd{\quad\hbox{and}\quad}

\def\v{{\mathcal V}}

\def\fm{(\cdot,\cdot)}
\def\a{\alpha}

\def\sub{\subseteq}

\def\lam{\lambda}

\def\1k{\frac{1}{k}}

\def\d{\delta}

\def\b{\beta}

\def\qed{\hfill$\Box$}

\def\sg{\sigma}

\def\hh{{\mathcal H}}

\def\bbbz{{\mathbb Z}}

\def\bbbf{{\mathbb F}}

\def\aa{\mathcal A}

\def\ll{\mathcal L}

\def\proof{\noindent {\bf Proof. }}

\def\jb{\bar{J}}
\def\cbox{\hfill$\lozenge$}

\def\ll{\mathcal L}
\def\hx{\hat{x}}
\def\hy{\hat{y}}
\def\hz{\hat{z}}
\def\bx{\bar{x}}
\def\by{\bar{y}}
\def\bz{\bar{z}}

\begin{document}

\author{S. Azam}
\thanks{The research was in part supported by a grant from IPM (No. 93160221).
This research is partially carried out in the IPM-Isfahan Branch.
The author would like to thank the Banach Algebra Center of Excellence for Mathematics,
University of Isfahan.}

\title[Kac--Moody--Malcev Superalgebras]{A New Characterization of Kac--Moody--Malcev Superalgebras}
\address
{School of Mathematics, Institute for
Research in Fundamental Sciences (IPM), P.O. Box: 19395-5746 and
Department of Mathematics\\ University of Isfahan\\Isfahan, Iran,
P.O. Box: 81745-163.} \email{azam@sci.ui.ac.ir; saeidazam@yahoo.com}

\date{}

\begin{abstract}
In the past two decades there has been a great attention to Lie (super)algebras which are extensions of affine Kac-Moody Lie (super)algebras, in certain typical or
axiomatic approaches. These Lie (super)algebras have been mostly studied under variations of the name "extended affine Lie (super)algebras".
We show that certain classes of Malcev (super)algebras  also can be put in this framework. This
in particular allows to provide new examples of Malcev (super)algebras which extend the known Kac-Moody Malcev (super)algebras.
\end{abstract}

\maketitle

\setcounter{section}{-1}
\section{\bf Introduction}\setcounter{equation}{0}
In 1990, two mathematical physicists,
H\o egh-Krohn and Torr\'esani \cite{hkt},
introduced a class  of Lie algebras in an axiomatic approach.
This new class generalizes the class of finite dimensional irreducible simple Lie algebras and affine Kac-Moody Lie algebras. Since then, this class and its certain extensions have been under intensive investigation, mostly, under the variations of the name extended affine Lie algebras (for a detailed account see \cite{aabgp} and \cite{neh}). Recently, this axiomatic approach has been extended to a superversion by M. Yousofzadeh \cite{you1}, introducing the class of extended affine Lie superalgebras, which covers extended affine Lie algebras, basic classical Lie superalgebras as well as affine Kac-Moody Lie superalgebras.

In the Lie algebra case, an extended affine Lie algebra over a filed $\bbbf$ of characteristic zero consists of three main ingredients, a Lie algebra $\ll$, an abelian subalgebra $\hh$ which with respect to $\ll$ has a root space decomposition with root system say $R$, and a symmetric non-degenerate invariant bilinear form $\fm$ whose restriction to $\hh$ is also non-degenerate. These three ingredients are related via two axioms; for each non-zero root  $\a$, there exist $\pm\a$-root vectors $x_{\pm\a}$ with $0\not=x_\a x_{-\a}\in\hh$, and that; the root vectors which are non-isotropic with respect to the form act locally nilpotently. The axioms for an extended affine Lie superalgebra are a superversion of these axioms.

The goal of  this note is two folded. First, to show that  by considering the above axioms in a general setting of algebras (see Definition \ref{main}), not necessarily Lie (super)algebras, one can treat more classes of algebras
including certain subclasses of Malcev (super)algebras  in
this setting. In this framework, we produce new examples of mostly infinite dimensional non-Lie Malcev (super)algebras which generalize the class of affine Kac-Moody Malcev (super)algebras.
Our second goal is to give a new prospective to the axioms of extended affine Lie (super)algebras. People familiar with the description of root systems of extended affine Lie (super)algebras (see \cite{aabgp}, \cite{you2}) know that these root systems enjoy certain  built in "finiteness" condition, namely the form restricted to the root system takes only a finite number of values.  We recognize this from the beginning in our axioms   (see axiom (E2) in Definition \ref{main}). This setting in particular allows to treat uniformly more interesting algebras related to the theory.

\section{\bf Set up}\setcounter{equation}{0}
Let $B$ be a finite integral domain. For $b\in B$, $\circ(b)$ denotes the order of $b$ in the abelian group $(B,+)$. Let $\aa=\oplus_{b\in B}\aa_b$ be a $B$-graded algebra. For a homogeneous element $x$ we denote its degree by $|x|$. Throughout this work, all elements of graded spaces are taken to be homogeneous. The algebra $\aa$ is called $B$-commutative if $xy=(-1)^{\circ(|x||y|)}yx$ for all $x,y$. A bilinear form
$\fm:\aa\times \aa\rightarrow\bbbf$ is called {\it invariant} if $(xy,z)=(x,yz)$ for all $x,y,z$, it is called {\it $B$-symmetric} if $(x,y)=-(-1)^{\circ(|x||y|)}(y,x)$ for all $x,y$, and is called {\it $B$-graded} if $(\aa_b,\aa_{b'})=\{0\}$ unless $b+b'=0$. If $B=\bbbz_2$, the terms $B$-commutative, $B$-symmetric and $B$-graded, will be called {\it supercommutative}, {\it supersymmetric} and {\it even}, respectively.
We denote the dual space of a vector space $\v$ by $\v^*$. To indicate that a remark or example is concluded we use the symbol $\lozenge$.

From now on we assume that $\aa$ is a $B$-commutative graded algebra with homogeneous spaces $\aa_{b}$, $b\in B$.

\begin{DEF}\label{deftoral}\em
Let $\aa$ be an algebra and $\hh$ a nonzero subalgebra of $\aa_0$.
For $\a\in\hh^*$ set $\aa^\a:=\{x\in\aa\mid hx=\a(h)x\hbox{ for all }h\in\hh\}.
$
We say that $(\aa,\hh)$ is a {\it toral pair} if
$\aa=\bigoplus_{\a\in\hh^*}\aa^\a$ and $\aa^0\aa^0\sub\aa^0$. In this case $\hh$ is called a {\it toral subalgebra}.
The set $R:=\{\a\in\hh^*\mid \aa^\a\not=0\}$ is called the {\it root system} of $(\aa,\hh)$.
Note that since $\hh\subseteq\aa_0$, for each $\a\in\hh^*$ we have $\aa^\a$ is $B$-graded, as a vector space, with $\aa^\a_b=\aa_\a\cap\aa^b,$ $b\in B$. For $b\in B$, we set $R_b:=\{\a\in R\mid \aa^\a_b\not=\{0\}\}$.
A toral pair $(\aa,\hh)$ is called {\it $B$-quadratic} if it is equipped with a
a $B$-symmetric, $B$-graded invariant form.
\end{DEF}

The following fact is worth mentioning in our general setting.

\begin{lem}\label{lem1}
Toral subalgebras are abelian, in the sense that all products are trivial.
\end{lem}

\proof Let $\hh$ be a toral subalgebra of $\aa$. Since $\hh\subseteq\aa_0$, $hh'=-h'h$ for all $h,h'\in\hh$, or equivalently $h^2=0$ for all $h$.
We are done if we show that for $h=\sum_\a h_\a$, $h_\a\in\aa^\a$, $h_\a=0$ for all $\a\not=0$.  Let $h'=\sum_\a h'_\a\in\hh$ be arbitrary.  We have
$$hh'=\sum_\a\a(h)h'_\a=-\sum_\a\a(h')h_\a.$$
Thus $\a(h)h'_\a=-\a(h')h_\a$ for all $\a$.
This forces $\a(h')h_\a=0$ if $h_\a\not=0$. Since $h'$ is arbitrary, we get $\a=0$ if $h_\a\not=0$.
\qed

\begin{DEF}\label{main}\em
Let $(\aa,\hh)$ be a $B$-quadratic toral pair with corresponding root system $R$ and bilinear form $\fm$. Then we call $(\aa,\hh)$ a $B$-{\it extended algebra} if the following two axioms hold:

(E1) for $b\in B$ and $\a\in R_b\setminus\{0\}$, there exists a pair $x_{\pm\a}\in\aa^{\pm\a}_{\pm b}$ with $0\not=x_\a x_{-\a}\in\hh$,

(E2) $\b(x_\a x_{-\a})$ varies in a  finite set for all $\a,\b\in R$ and all pairs $x_{\pm\a}\in\aa^{\pm\a}$ with $(x_\a,x_{-\a})=1$.

\noindent The $B$-extended algebra $(\aa,\hh)$ is called a {\it $B$-extended affine algebra} if

(E3) the form $\fm$ is nondegenerate on $\aa$ and on $\hh$.

If ${\mathcal P}$ is a homogenous variety of algebras, then a  $B$-extended (affine) algebra $(\aa,\hh)$ is called a {\it $B$-extended ${\mathcal P}$-algebra} if $\aa$ belongs to ${\mathcal P}$. If $B$ is the trivial grading or $B$ is understood from the context, then we drop $B$ form the terminology. For example, with this terminology, a Lie (super)algebra $\aa$ satisfying (E1)-(E3) is called an extended affine Lie (super)algebra.
\end{DEF}

{\bf Notation.}  Following \cite[Definition 1.3]{stu}, we call a pair $x_{\pm\a}$ as in (E1) a {\it test-pair}.
 We also call a test-pair $x_{\pm\a}$ a {\it block test-pair} if $(x_\a,x_{-\a})=1$. In extended affine Lie (super)algebras block test-pairs are responsible for the existence of $\mathfrak{sl}_2(\bbbf)$-blocks in the Lie case and $\mathfrak{osp}(1,2)$-blocks in Lie super case. The main structural features of the cores of Lie (super)algebras are essentially deduced from the block test-pairs. With these terminologies, extended algebras are quadratic toral pairs for which corresponding to each non-zero root there exists a test-pair, and roots take finite values on all block test-pairs.


For people familiar with the class of affine Lie (super)algebras and its generalizations, axioms (E1) and (E3) looks quite natural. In part (ii) of the remark \ref{rrem1} below, we discuss axiom (E2) and show that when the form on $\aa$ restricted to $\hh$ is non-degenerate, this axiom can be replaced with a more natural and simpler expression (see Remark \ref{rrem1}(ii)).

\begin{rem}\label{rrem1}\em Let $(\aa,\hh)$ be  a $B$-quadratic extended algebra with root system $R$ and bilinear form $\fm$.
Consider the map
$$\nu:\hh\longrightarrow\hh^*,\quad h\mapsto (h,\cdot).
$$
If $\nu$ is one to one, we transfer the form on $\hh$ to $\nu(\hh)$ by
\begin{equation}\label{form1}
(\a,\b):=(t_\a,t_\b)\hbox{ where }t_\a:=\nu^{-1}(\a),\quad\a\in\nu(\hh).
\end{equation}
We have now the following four remarks.

(i) By definition, the form $\fm$ is symmetric on $\aa_0$ and so on $\hh$. Assume the form $\fm$ restricted to $\hh$ is nondegenerate, then the map $\nu$ is injective and so we can transfer it to $\nu(\hh)$ by (\ref{form1}).
Now for any test-pair $x_{\pm\a}$, $\a\in\nu(\hh)$, and $h\in\hh$ we have,
$$(h,x_\a x_{-\a})=(hx_\a,x_{-\a})=\a(h)(x_\a,x_{-\a}).$$
Therefore, as the $\fm$  is nondegenerate on $\hh$, we get $x_\a x_{-\a}=(x_\a,x_{-\a})t_\a$. In particular, if $\a\not=0$, then $(x_\a, x_{-\a})\not=0$. Therefore if (E1) holds,
 the non-degeneracy of the form on $\hh$ implies that
$R\sub\nu(\hh)$.

(ii)  Assume the form on $\hh$ is non-degenerate. Let $\a,\b\in R$ and $x_{\pm\b}$ be a block test-pair. We have from part (i) that $\a(x_\b x_{-\b})=\a(t_\b)=(\a,\b)$. Therefore, under the nondegeneracy of the form on $\hh$,  the axiom (E2) is equivalent to

{(E2)$'$ the set $\{(\a,\b)\mid\a,\b\in R\}$ is finite.}

\noindent
In particular, in the definition of a $B$-extended affine algebra, (E2) can be replaced with (E2)$'$.

(iii) Assume $(\aa,\hh)$ satisfies (E1)-(E3). Since the form is invariant and $B$-graded, and  $hx=-hx$ for $h\in\hh,\;x\in\aa$, it follows that $(\aa^\a_\sg,\aa^\b_\tau)=\{0\}$ unless $\a+\b=0$ and $b+b'=0$. In particular, the form restricted to $\aa^\a_b+\aa^{-\a}_{-b}$, $\a\in\hh^*$, $b\in B$, is nondegenerate, and $R_b=-R_b$ for each $b\in B$.

(iv) People familiar with affine Lie algebras and their extensions, such as extended affine Lie (super)algebras, realize that the axiom (E2) (or (E2$'$)), is a replacement for the usual ad-nilpotency axiom which states that for $\a\in R$ with $(\a,\a)\not=0$, the elements in $\aa^\a$ are locally nilpotent. In the presence of the Jacobi (super)identity the algebra $\aa$ is $\hh^*$-graded and so this axiom implies further consequences such as boundedness of root strings. However, one knows that in general non-associative algebras (toral pairs) such as Malcev (super)algebras,  the algebra $\aa$ is not in general $\hh^*$-graded (see Remark \ref{rem33}(ii)). This gives a justification that our axiom (E2) (or (E2$'$)) is more natural in the context of general non-associative algebras.\cbox
\end{rem}

\begin{ex}\label{rexa1}\em (i) One knows that simple finite dimensional Lie algebras, affine Lie algebras and their various generalizations such as extended affine Lie (super)algebras as well as invariant affine reflection algebras are examples of $B$-extended affine algebras, where $B$ is the trivial grading in the Lie algebra case and the $\bbbz_2$-grading in the Lie super case.

(ii) Locally finite split Lie algebras (\cite{NS}) with integrable roots are examples
of extended affine algebras.

(iv) Let $G$ be a torsion free abelian group and $\aa$ be a Lie $G$-torus in the sense of \cite{yos1} of \cite{yos2}.
By \cite[Theorem 7.1]{yos1}, $\aa$ admits a non-zero symmetric invariant form.
From this, the involved irreducible finite root system, and the division property it follows that a Lie $G$-torus is an extended Lie algebra. The same is through if one replace the involved finite root system with an irreducible locally finite
root system.\cbox

\end{ex}

\begin{DEF}\label{def14}\em
For a $B$-quadratic toral pair $(\aa,\hh)$ denote by $\aa_c$ the subalgebra generated by root spaces $\aa^\a$, $\a\in R^\times$ where
$$R^\times:=\{\a\in R\mid \a (x_\b x_{-\b})\not=0\hbox{ for some test-pair } x_{\pm\b}\in\aa^{\pm\b}\}$$ is the set of {\it non-isotropic} roots of $R$. If the form $\fm$ is non-degenerate on $\hh$ then as Remark \ref{rrem1} shows, we have $R^\times=\{\a\in R\mid (\a,R)\not=\{0\}\}.$ We call $\aa_c$, the {\it core} of $\aa$, and we set $\hh_c:=\hh\cap\aa_c.$
\end{DEF}

\begin{pro}\label{yy1}
Let $(\aa,\hh)$ be an extended affine Lie superalgebra, with Cartan subalgebra $\hh$ and root system $R$.
Then $(\aa_c,\hh_c)$ is an extended superalgebra.
\end{pro}

\proof
One can easily see that $\aa_c$ is an ideal of $\aa$ with $\aa_c=\sum_{\a\in\hh^*}\aa_c\cap\aa^\a$. Let $\hh_c:=\hh\cap\aa_c$ and for $\a\in (\hh_c)^*$, let $\aa^\a_c:=\{x\in\aa_c\mid hx=\a(h)x\hbox{ for all }h\in\hh_c\}$. Then $\aa_c=\sum_{\a\in(\hh_c)^*}\aa_c^\a$ and for $\a\in(\hh_c)^*$, $\aa_c^\a=\sum\aa^\b\cap\aa_c$ where sum runs over all $\b\in R$ with $\b=\a$  restricted to $\hh_c$. Moreover $\aa^0_c\aa^0_c\sub\aa^0_c$.
Thus $(\aa_c,\hh_c)$ is a quadratic toral pair. Let $R_c$ be the root system of $R$ and $0\not=\a\in (R_c)_b$, $b\in B=\bbbz_2$.  Then there exists $\b\in R$ such that $\b=\a$ on $\hh_c$ and $\{0\}\not=\aa^\b\cap\aa_c\sub\aa_c\cap\aa^\a\cap\aa_b$. Then $0\not=\b\in R_b$. We have
$$\hh_c=\hh\cap\aa_c\sub\aa^0\cap\aa_c=\sum_{\gamma\in R^\times}\aa^\gamma\aa^{-\gamma}\sub\sum_{\gamma\in R^\times}\bbbf t_\gamma.
$$
Therefore, since $\b(\hh_c)=\a(\hh_c)\not=\{0\}$, there exists $\gamma\in R$ with $\b(t_\gamma)=(\b,\gamma)\not=0$. In particular $\b\in R^\times$.

Since $0\not=\b\in R_b$, we have from (E1) that there exists a test-pair $x_{\pm\b}\in\aa^{\pm\b}_b$. Since  $\b\in R^\times$, we have $x_{\pm\b}\in\aa^{\pm\b}_c$ and $x_\b x_{-\b}\in\hh_c$. Thus
(E1) holds for $(\aa_c,\hh_c)$.

Next, we show that (E2) holds for $(\aa_c,\hh_c)$. Let $\a\in R_c$ and $x_{\pm\a}$ be a block test-pair. Then as we have seen,
$x_\a=\sum_{\b}x^\a_\b$, $x_{-\a}=\sum_{\b}x^{-\a}_{-\b}$ where $\b$ runs over the subset $R_\a$ of $R$ consisting of roots whose restriction to $\hh_c$ coincide with $\a$. Since $x_\a x_{-\a}\in\hh_c\subseteq\aa^0$, and the form is $\hh^*$-graded, we have
$x_\a x_{-\a}=\sum_{\b\in R_\a} x_\b^\a x^{-\a}_{-\b}$ and $1=(x_\a,x_{-\a})=\sum_{\b\in R_\a}(x_\b^\a,x_{-\b}^{-\a}).$ Now let $\gamma\in R_c$. Since $\aa^\gamma_c\not=\{0\}$, there exists $\tau\in R$ which restricted to $\hh_c$ coincides with $\gamma$.
Then
$$\gamma(x_\a x_{-\a})=\sum_{\b\in R_\a}\tau(x_\b^\a x_{-\b}^{-\a})=\sum_{\b\in R_\a}(x_\b^\a, x_{-\b}^{-\a})\tau(t_\b)=\sum_{\b\in R_\a}(x_{\b}^\a, x_{-\b}^{-\a})(\tau,\b).$$
Now from this and the facts that $\sum_{\b\in R_\a}(x_{\b}^\a, x_{-\b}^{-\a})=1$ and $\{(\tau,\b)\mid\b\in R\}$ is finite, it follows that (E2) holds for $(\aa_c,\hh_c)$. Thus  $(\aa_c,\hh_c)$ is an extended superalgebra.\qed

\section{\bf Affine construction}\setcounter{equation}{0}\label{affine-construction}
Let $\aa$ be a $B$-commutative algebra equipped with a $B$-symmetric $B$-graded invariant bilinear form. Let $G$ be a non-trivial
torsion  free abelian group and $\lam:G\times G\rightarrow\bbbf^\times$ a $2$-cocycle satisfying $\lam(\sg,\tau)=\lam(\tau,\sg)$ for all $\sg$, $\tau$,  and $\lam(0,0)=1$. Let $\bbbf^t[G]$ be the commutative associative torus associated with the pair $(G,\lam)$ (see \cite[\S 1.4]{ayy}). One knows that as a vector space $\bbbf^t[G]$ has a $\bbbf$-basis $\{t^\sg\mid\sg\in G\}$ with algebra multiplication given
by $t^\sg t^\tau=\lam(\sg,\tau)t^{\sg+\tau}$, $\sg,\tau\in G$.
Consider the tensor algebra
$$\ll({\aa}):=\aa\otimes\bbbf^t[G].$$
Then $\ll(\aa)$ is a
$B$-graded algebra with $\ll(\aa)_b=\aa_b\otimes\bbbf^t[G]$ for $b\in B$.

The form on $\aa$ can be extended to a form on $\ll(\aa)$ by
$$(x\otimes t^\sg,y\otimes t^\tau):=(x,y)\lam(\sg,\tau)\d_{\sg,-\tau}.$$
It can be easily checked that  $\ll(\aa)$ is a $B$-commutative algebra, the extended form is $B$-symmetric on $\ll(\aa)$, and that
the form $\fm$ is non-degenerate on $\ll(\aa)$ if it is non-degenerated on $\aa$.

Identify $G$ with $G\otimes 1$ in $G\otimes_\bbbz\bbbf$. Let
$$\v:=G\otimes_{\bbbz}\bbbf.$$  We may consider $\v^*$ as a subspace of
derivation algebra $\hbox{Der}\ll(\aa)$ of $\ll(\aa)$ by
$$d(x\otimes t^\sg)=d(\sg)x\otimes t^\sg,\quad (d\in\v^*,\;\sg\in G).
$$
Let
$$\tilde{\ll}(\aa):=\ll(\aa)\oplus\v\andd \hat{\ll}(\aa):=\aa\otimes\bbbf^t[G]\oplus\v\oplus\v^*.$$ Then both these algebras are $B$-graded with
$$\begin{array}{c}\tilde{\ll}(\aa)_0=(\aa_0\otimes\bbbf^t[G])\oplus\v\andd\tilde{\ll}(\aa)_b=\aa_b\otimes\bbbf^t[G],\hbox{ for }b\not=0,\\
\hat{\ll}(\aa)_0=(\aa_0\otimes\bbbf^t[G])\oplus\v\oplus\v^*\andd \hat{\ll}(\aa)_b=\aa_b\oplus\bbbf^t[G],\hbox{ for }b\not=0.\end{array}$$
We then turn $\tilde{\ll}(\aa)$ and $\hat{\ll}(\aa)$ into algebras by assuming that $\v$ is central, $\v^*\v^*=\{0\},$ and that for $x,y\in\aa$, $\sg,\tau\in G$ and $d\in\v^*$,
$$\begin{array}{c}
(x\otimes t^\sg)(y\otimes t^\tau):=\lam(\sg,\tau)\big(xy\otimes t^{\sg +\tau}+\d_{\sg,-\tau}(x,y)\sg\big),\\
d\cdot (x\otimes t^\sg)=-(x\otimes t^\sg)\cdot d:=d(x\otimes t^\sg).
\end{array}$$

The form on $\ll(\aa)$ can be extended to a  bilinear form on both $\tilde{\ll}(\aa)$ and $\hat{\ll}(\aa)$ by $(\v,\tilde{\ll}(\aa)\oplus\v)=(\v^*,\tilde{\ll}(\aa)\oplus\v^*)=\{0\}$ and $(d,\sg)=d(\sg)$ for $d\in\v^*$ and $\sg\in\v$.

For $X=x\otimes t^\sg$, $Y=y\otimes t^\tau$, $x,y\in\aa$ and $\sg,\tau\in G$, we have
\begin{eqnarray*}
XY&=&\lam(\sg,
\tau)\big(xy\otimes t^\sg t^\tau+(x,y)\d_{\sg,-\tau}\sg)\big)\\
&=&\lam(\sg,\tau)\big((-1)^{\circ(|x||y|)}yx
\otimes t^\tau t^\sg+(-1)^{\circ(|x||y|)}(y,x)\d_{\sg,-\tau}\tau\big)\\
&=&(-1)^{\circ(|X||Y|)}YX,
\end{eqnarray*}
and
$$(X,Y)=\lam(\sg,\tau)(x,y)\d_{\sg,-\tau}=-(-1)^{\circ(|x||y|)}\lam(\sg,\tau)(yx)=-(-1)^{\circ(|X||Y|)}(Y,X).$$
It follows that both $\tilde{\ll}(\aa)$ and $\hat{\ll}(\aa)$ are $B$-commutative and the extended form is $B$-symmetric.

Assume further that $\aa$ has a subalgebra $\hh$ such that $(\aa,\hh)$ is a toral pair with root system $R$.
Set $\tilde{\hh}:=(\hh\otimes 1)\oplus\v\sub\tilde{\ll}(\aa)_0$ and $\hat{\hh}:=(\hh\otimes 1)\oplus\v\oplus\v^*\sub\hat{\ll}(\aa)_0.$ We identify $\tilde{\hh}^*$ with $\hh^*\oplus\v^*$ and
$\hat\hh^*$ by $\hh^*\oplus\v^*\oplus(\v^*)^*$.
We also identify $\v$ as a subspace of $(\v^*)^*$ in the natural way.
For $\hat\a\in\hat\hh^*$ define $\hat\ll(\aa)^\a$ in the usual manner. Then we get $\hat\ll(\aa)=\sum_{\hat\a\in\hat R}\hat\ll(\aa)^{\hat\a}$ with $\hat{R}=R+G$ and for $\a\in R$, $\sg\in G$, $\hat\ll(\aa)^{\a+\sg}=\aa^\a\otimes t^\sg$ if $\a+\sg\not=0$ and
$\hat{\ll}(\aa)^0=\aa^0\otimes 1\oplus\v\oplus\v^*.$
We also have $\hat\ll(\aa)^0\hat\ll(\aa)^0\subseteq\aa^0\aa^0\otimes 1\sub\aa^0\otimes 1\sub\hat\ll(\aa)^0$. Thus $(\hat{\ll}(\aa),\hat\hh)$ is a quadratic toral pair. Similarly, we see that $(\tilde\aa,\tilde\hh)$ is a $B$-quadratic  toral pair with root system $R$, where $\tilde{\ll}(\aa)^\a=\aa^\a\otimes \bbbf^t[G]$ if $0\not=\a\in R$ and $\tilde{\ll}(\aa)^0=\aa^0\otimes\bbbf^t[G]\oplus\v$.

The following result extends the  realization of untwisted affine  Kac-Moody Lie algebras to a much broader setting (compare with \cite[\S 7]{ahy}).  

\begin{pro}\label{nnn1} Consider the $B$-commutative algebra $\aa$ and assume that $(\aa,\hh)$ is a $B$-quadratic toral pair with root system $R$ and the bilinear form $\fm$. Let $G$ be a non-trivial torsion free abelian group and $\lam:G\times G\rightarrow \bbbf^{\times}$ a symmetric quadratic map.

(i) The pairs $(\tilde{\ll}(\aa),\tilde{\hh})$ and $(\hat{\ll}(\aa),\hat{\hh})$ are quadratic toral pairs with root systems $R$ and $\hat{R}=R+G$, respectively. Moreover, if the form on $\aa$, respectively $\hh$ is non-degenerate, then so is the extended form on $\hat{\ll}(\aa)$, respectively $\hat\hh$.

(ii) Suppose $(\aa,\hh)$ is a $B$-extended algebra. Then

$\quad$(a)  $(\tilde{\ll}(\aa),\tilde\hh)$ is a $B$-extended algebra, in particular if $(\aa,\hh)$ is an extended affine algebra then so is $({\ll}(\aa),\hh)$.

$\quad$(b) If $\aa^{0}_{0}=\hh$ and form is non-degenerate on $\aa^{0}$, then  $(\hat{\ll}(\aa),\hat{\hh})$ is a $B$-extended algebra. In this case if $(\aa,\hh)$ is a $B$ extended affine algebra, then so is $(\hat{\ll}(\aa),\hat\hh)$.

\end{pro}

\proof
(i) The details is given in the paragraph preceding the proposition.

(ii) Assume that (E1)-(E2) hold for $(\aa,\hh)$, $\aa^0_0=\hh$ and that the form restricted to $\aa^0$ is non-degenerate. We show that (E1)-(E2) hold for $(\hat{\ll}(\aa),\hat{\hh})$. Let $b\in B$ and $\hat{\a}\in\hat{R}_b\setminus\{0\}$. Since $\hat{R}_b=R_b+G$, we have $\hat{\a}=\a+\sg$ for some $\a\in R_b$ and $\sg\in G$. We divide the argument in two cases $\a=0$ and $\a\not=0$. Assume first that $\a=0$ (forcing $\sg\not=0$). Since $\a=0\in R_b$, $\aa^0_b\not=0$. Since the form in non-degenerate on $\aa^0$ and $B$-graded, there exist $x_\pm\in\aa^0_{\pm b}$ with $(x_+,x_-)\not=0$.
Now $x_\pm\otimes t^{\pm\sg}\in\hat{\ll}(\aa)^{\pm\sg}_{\pm b}=\aa^0_{\pm b}\otimes t^{\pm\sg}$. Then since $\aa^0\aa^0\sub\aa^0$, we have
\begin{eqnarray*}
X:=(x_+\otimes t^\sg)(x_-\otimes t^{-\sg})&=&\lam(\sg,-\sg)x_+x_-\otimes 1+
\lam(\sg,-\sg)(x_+,x_-)\sg\\
&\in&(\aa^0\aa^0\cap\aa_0)\otimes 1\oplus\v\\
&\sub&(\aa^0\cap\aa_0)\otimes 1\oplus\v\\
&=&\hh\otimes 1\oplus\v\sub\hat{\hh}.
\end{eqnarray*}
But $X\not=0$ as $(x_+,x_-)\not=0$ and  $d_i(\sg)\not=0$ for at least one $(i,\sg)\in I\times G$. Next, assume that $\a\not=0$. By (EA1), there exists $x_\pm\in\aa^{\pm\a}_{\pm\b}$ such that $0\not=x_+x_-\in\hh$. Then $x_\pm\otimes t^{\pm\sg}\in\hat{\ll}(\aa)^{\pm\hat{\a}}_{\pm b}$ and
$$(x_+\otimes t^\sg)(x_-\otimes t^{-\sg})= \lam(\sg,-\sg)x_+x_-\otimes 1+\lam(\sg,-\sg)(x_{-},x_{+})\sg,$$ which is a non-zero element of $\hat{\hh}$ as required. Thus (E1) holds for $(\hat{\ll}(\aa),\hat\hh)$.

Next, we show that (E2) holds. Let $\hat\a,\hat\b\in\hat R$, $b\in B$, and $\hat{x}_{\pm\hat\a}\in\hat{\ll}(\aa)^{\pm\hat\a}_{\pm b}$ be a block test-pair.
Since $\hat{\ll}(\aa)^0_0=\aa^0_0\oplus\v\oplus\v^*=\hh\oplus\v\oplus\v^*$, we have $\hat{x}_\a\hat{x}_{-\a}=0$.
So to check (E2), we may assume that $(\hat\a,b)\not=(0,0)$. In this case $\hat{x}_{\pm\hat{\a}}=x_{\pm\a}\otimes t^{\pm\sg}$ for some
$\a\in R$, $\sg\in \v$. Now let $\hat\b=\b+\tau\in\hat R$, with $\b\in R$, $\tau\in\v\sub(\v^*)^*$. Then
$$\hat\b(\hat{x}_\a\hat{x}_{-\a})=(\b+\tau)(x_\a x_{-\a}\otimes 1)\lam(\sg,-\sg)=\b(x_\a x_{-\a})\lam(\sg,-\sg).$$
Since $1=(\hat x_\a,\hat{x}_{-\a})=(x_\a,x_{-\a})\lam(\sg,-\sg)$ and (E2) holds for $(\aa,\hh)$, we see that (E2) holds for $\hat{\ll}(\aa)$, that is $(\hat{\ll}(\aa),\hat\hh)$ is a $B$-extended algebra.

Next, we show that $(\tilde{\ll}(\aa),\tilde\hh)$ satisfies (E1)-(E2).
Assume first that $0\not=\a\in R_b$, $b\in B$. Since (E1) holds for $\aa$, there exists a test-pair $x_\pm\in\aa^\pm_b$. Then $x_\pm\otimes 1$ is a test-pair in $\tilde{\ll}(\aa)^\pm_b$, and so (E1) holds. Next let $\a\in R$ and that $\hat{x}_\pm\in\tilde{\ll}(\aa)^{\pm\a}$ be a block test-pair. Then $\hat{x}_+=\sum_ix^i_\a\otimes t^{\sg_i}$ (mod $\v$) and
$\hat{x}_+=\sum_ix^i_{-\a}\otimes t^{\tau_i}$ (mod $\v$), for some $x^i_{\pm\a}\in\aa^{\pm\a}$ and $\sg_i,\tau_i\in G$.
Since $\sum_{i,j}x^i_\a x^j_{-\a}\otimes t^{\sg_i}t^{\tau_j}\in\hh\otimes 1$ (mod $\v$), we may assume $\hat{x}_{\pm}=\sum_i x^i_{\pm\a}\otimes t^{\pm\sg_i}$ (mod $\v$) where $x_{\pm\a}^i\in\aa^{\pm\a}$, $\sg_i$'s are distinct elements of $G$ and $\sum_{i\not= j}x^i_\a x^j_{-\a}=0$ .
Set $X_{\pm}:=\sum_{i}x^i_{\pm\a}\in\aa^{\pm\a}$. Then
$0\not=X_+X_-=\sum_{i}x_\a^ix_{-\a}^i\in\hh$ and so the pair $X_\pm$ is a test-pair in $(\aa,\hh)$. Moreover,
$$(X_+,X_-)=\sum_{i}(x^i_\a,x^i_{-\a})=\sum_{i,j}(x^i_\a,x^j_{-\a})\d_{\sg_i,-\sg_j}=(\hat{x}_+,\hat{x}_-)=1.$$
Now for any $\b\in R$, $\b(\hat{x}_+\hat{x}_{-})=\b(X_+X_-\otimes 1 +\sg)=\b(X_+X_-)$ for some $\sg\in\v$. Thus (E2) holds for $\tilde{\ll}(\aa)$ as it holds for $\aa$.

Finally, as we have already seen if the form on $\aa$ is non-degenerate, then is so  on $\ll(\aa)$, $\hat{\ll}(\aa)$. Similarly if the form on $\hh$ is non-degenerate, then it  is also non-degenerate on $\hat{\hh}$. This completes the proof of (ii).\qed

\section{\bf Extended Malcev (super)algebras}\setcounter{equation}{0}
In this section, we first review some basic facts about Malcev (super)algebras and then use our affine construction given is Section \ref{affine-construction} to produce new examples of extended  non-Lie Malcev (super)algebras. Our method in part gives a new characterization to the known infinite dimensional super-Kac-Moody Malcev algebras (see \cite{CRT} and \cite{O}).

The Malcev algebras (or Moufang-Lie algebras) were first introduced by A. I. Malcev in 1955 [Mal]  as an algebra satisfying identities
$$x^2=0$$ and
\begin{equation}\label{malcev1}
J(x,y,xz)=J(x,y,z)x\end{equation}
for all x,y,z in $\aa$, where
$$J(x,y,z):=(xy)z+(zx)y+(yz)x$$
is the Jacobian of $x,y,z$.
The identity (\ref{malcev1}) is called the {\it Malcev identity} which can be written in the form
\begin{equation}\label{malcev2}
((xy)z)x+((yz)x)x+((zx)x)y=(xy)(xz).
\end{equation}
If the characteristic of the field is not equal to $2$, which is in our case, one can see after some linearizations  that in an anticommutative algebra the Malcev identity is
equivalent to the Sagle identity
$$((xy)z)t+((tx)y)z+((zt)x)y+((yz)t)x=(xz)(yt).
$$
The Sagle identity is linear in each variable and is invariant under cyclic permutation of its variables.

Clearly any Lie algebra is a Malcev algebra.
A {\it binary Lie algebra} is an algebra such that any of two elements of which generate a Lie subalgebra. A binary Lie algebra is given by the system of identities
$$x^2=J(xy, x,y)=0.$$
Therefore, any Malcev algebra is a binary Lie algebra as
$$J(xy,x,y)=J(y,x,yx)=J(y,x,x)y=0.$$
The class of Malcev algebras can be regarded as a class between Lie algebras and binary Lie algebras. Any alternative algebra $A$
turns into a Malcev algebra $A^-$ with the commutator product as its algebra structure. The algebra $A^-$ of a Cayley-Dickson algebra $A$ is a
Malcev algebra which is not a Lie algebra. While this is known by PBW theorem that any Lie algebra is a subalgebra of $A^-$ for some associative algebra,
it is an open problem whether or not any
Malcev algebra is a subalgebra of $A^-$ for some alternative algebra $A$.

In the same way that Lie algebras appear as algebraic models for Lie groups, Malcev algebras are the algebraic structures associated with analytic ``{\it Moufang loops}''.
More precisely, the tangent space $T_e(G)$ of a Moufang loop $G$ at identity, equipped with anti-commutative multiplication,
is a Malcev algebra.

For a $\bbbz_2$-graded  superalgebra $\aa$, the superjacobian is defined by
$$\jb (x,y,z):=(-1)^{|x||z|}(xy)z+(-1)^{|z||y|}(zx)y+(-1)^{|y||x|}(yz)x.$$

Then a Malcev superalgebras is a superalgebra satisfying the superversion of the identities (\ref{malcev1}) for homogenous elements, namely
$$xy=-(-1)^{|x||y|}yx\andd \jb (x,y,xz)=\jb (x,y,z)x.$$
In terms of a superversion of the Sagle identity, a Malcev superalgebra is a superalgebra satisfying
\begin{equation}\label{malcev3}
\begin{array}{c}
xy=-(-1)^{|x||y|}yx,\vspace{2mm} \\
(-1)^{|y||z|}(xz)(yt)=((xy)z)t+(-1)^{|x|(|y|+|z|+|t|)}((yz)t)x\\\qquad\qquad+(-1)^{(|x|+|y|)(|z|+|t|)}((zt)x)y+(-1)^{|t|(|x|+|y|+|z|)}((tx)y)z
\end{array}
\end{equation}
for homogenous elements.

\begin{rem}\label{rem33}\em
(i) The term ``extended algebras'' is used in the literature in another contexts too.  In fact an algebra satisfying the identities
\begin{equation}\label{malcev4}xy=-yx\andd J(x,y,xy)=0,
\end{equation}
is called an {\it extended algebra} by A. A. Sagle \cite{sa2}. Lie algebras and Malcev algebras are subclasses of extended algebras. Over an algebraically closed field of characteristic zero, a simple finite dimensional  extended algebra (in the sense of Sagle)  is a simple Lie algebra or the simple seven dimensional Malcev algebra if and only if the trace form, $(x,y) := \hbox{trace}RxRy$ is a non-degenerate invariant form, where $R_x$ is the right multiplication by $x$ (see \cite{sa2}).

(ii) Let $(\aa,\hh)$ be a toral pair, where $\aa$ satisfies the identities (\ref{malcev4}). Then as in \cite{sa2}, one can prove that
\begin{equation}
\aa^\a \aa^\b\sub \aa^{\a+\b}\hbox{ if }\a\not=\b.
\end{equation}
Moreover, if for each $\a\not=0$ we have $\aa_\a\aa_\a\sub\sum_{\b\not=0}\aa_\b$, for example if $\aa$ is a Malcev algebra, then we have
$$\aa^\a\aa^\a\sub \aa^{2\a}+\aa^{-\a}.$$\cbox
\end{rem}

\begin{ex}\label{non-lie}\em
(A non-Lie example.) Consider the seven dimensional algebra $\aa$ with basis $e_1,\ldots,e_7$ and product  given by the table
$$\begin{array}{c|ccccccc}
\cdot&e_1& e_2 &e_3&e_4&e_5&e_6&e_7
\\
\hline
e_1&0&2e_2&2e_3&2e_4&-2e_5&-2e_6&-2e_7\\
\hline
e_2&-2e_2&0&2e_7&-2e_6&e_1&0&0\\
\hline
e_3&-2e_3&-2e_7&0&2e_5&0&e_1&0\\
\hline
e_4&-2e_4&2e_6&-2e_5&0&0&0&e_1\\
\hline
e_5&2e_5&-e_1&0&0&0&-2e_4&2e_3\\
\hline
e_6&2e_6&0&-e_1&0&2e_4&0&-2e_2\\
\hline
e_7&2e_7&0&0&-e_1&-2e_3&2e_2&0\\
\end{array}
$$
Taking $\hh:=\bbbf e_1$, we see that $(\aa,\hh)$ is an extended affine algebra with
$R=\{0,\pm\a\}\sub\hh^*$ where $\a(e_1)=2$. In fact $\aa^\a=\bbbf e_2\oplus\bbbf e_3\oplus\bbbf e_4$ and $\aa^{-\a}=\bbbf e_5\oplus\bbbf e_6\oplus\bbbf e_7$. The algebra $\aa$ is equipped with the invariant nondegenerate bilinear form given by
$$
\big((e_i,e_j)\big)=\left(\begin{array}{ccccccc}
1&0&0&0&0&0&0\\
0&0&0&0&1&0&0\\
0&0&0&0&0&1&0\\
0&0&0&0&0&0&1\\
0&1&0&0&0&0&0\\
0&0&1&0&0&0&0\\
0&0&0&1&0&0&0\\
\end{array}
\right).
$$
This is the basic example of a  non-Lie Malcev algebra obtained from minus algebra of trace zero elements in an octonion algebra
(see \cite[Example 3.2]{sa1}). We see here that, unlike the Lie case, we have $\aa^\a\aa^\a\sub\aa^{-\a}$. As this example shows Malcev toral pairs ($\aa,\hh)$ are not in general $\hh^*$-graded, see Remark \ref{rem33}(ii).\cbox
\end{ex}

Let $G$ be a non-trivial torsion free abelian group and $\lam:G\times G\rightarrow \bbbf^{\times}$ a symmetric quadratic map. According to Proposition \ref{nnn1}, if $(\aa,\hh)$ is a $B$-extended algebra  then
$(\tilde{\ll}(\aa),\tilde\hh)$ is also a $B$-extended algebra, in particular if $(\aa,\hh)$ is an extended affine algebra then so is $({\ll}(\aa),\hh)$.
Moreover, if $\aa^{0}_{0}=\hh$ and the form is non-degenerate on $\aa^{0}$, then  $(\hat{\ll}(\aa),\hat{\hh})$ is a $B$-extended algebra. In this case if $(\aa,\hh)$ is a $B$ extended affine algebra, then so is $(\hat{\ll}(\aa),\hat\hh)$. In the following proposition we investigate under which circumstances the Malcev superstructure of $\aa$ can be transferred to the affine constructions $\tilde{\ll}(\aa)$ and $\hat{\ll}(\aa)$.

\begin{pro}\label{nnn1} Let $\aa$ be a Malcev superalgebra. Let $G$ be a non-trivial torsion free abelian group and $\lam:G\times G\rightarrow \bbbf^{\times}$ a symmetric quadratic map.
Then $\tilde{\ll}(\aa)$ is a Malcev superalgebra. Moreover, $\hat{\ll}(\aa)$ is a Malcev superalgebra if and only if $\aa$ is a Lie superalgebra, in which case $\hat{\ll}(\aa)$ is a Lie superalgebra.
\end{pro}

\proof
 We first consider $\hat{\ll}(\aa)$. Let $\hat{x}:=x\otimes t^{\bar{x}}+v_x+d_x$, $\hat{y}:=x\otimes t^{\bar{y}}+v_y+d_y$ and $\hat{z}:=z\otimes t^{\bar z}+v_z+d_z$ be three arbitrary homogenous elements of $\hat{\ll}(\aa)$, where $x,y,z\in\aa$, $\bar x,\bar y, \bar z\in G$, $v_x,v_y,v_z\in\v=\bbbf\otimes G$ and $d_x,d_y,d_z\in\v^*$. We note that, if for example $v_x\not=0$, or $d_x\not=0$, then $|x|=0$, as $\hat{x}$ is homogeneous.
Then, using the facts that $\aa$ is supercommutative and the form is even,  a direct computation gives
$$
\bar{J}(\hat x,\hat y,\hat x\hat z)=A+B+C,
$$
with
$$A=\lam(\bar x,\bar y)\lam(\bar x,\bar z)\lam(\bar x+\bar y,\bar x+\bar z)\bar J(x,y,xz)\otimes t^{2\bar x+\bar y+\bar z}
$$
$$B=-\lam(\bx,\by)(\lam(\bx+\by,\bx)d_z(\bar x)\bar J(x,y,x)\otimes t^{2\bar x+\bar y},
$$
and
$$C=b\otimes t^{\bar x+\bar y+\bar z},$$
where
\begin{eqnarray*}
b&=&\alpha d_x(\bar z)(xy)z\lam(\bar x,\bar y)\lam(\bar x+\bar y,\bar z)\\
&&+\alpha d_x(\bar y)y(xz)\lam(\bar x,\bar z)\lam(\bar x+\bar z,\bar y)\\
&&-\b  d_x(\bar x+\bar z)(xz)y\lam(\bar x,\bar z)\lam(\bar x+\bar z,\bar y)\\
&&+\b  d_x(\bar z)(zx)y\lam(\bar x,\bar z)\lam(\bar x+\bar z,\bar y)\\
&&+\gamma  d_x(\bar z)(yz)x\lam(\bar y,\bar z)\lam(\bar y+\bar z,\bar x)\\
&&-\gamma  d_x(\bar x+\bar y+\bar z)y(xz)\lam(\bar x,\bar z)\lam(\bar x+\bar z,\bar y),
\end{eqnarray*}
and
$$\a=(-1)^{|x||xz|},\quad\b=(-1)^{|y||xz|}\andd \gamma=(-1)^{|x||y|}.
$$

Similarly, we get
$$J(\hat x,\hat y,\hat z)\hat x=A'+B'+C',
$$
with
$$A'=\lam(\bar x,\bar y)\lam(\bar x+\bar y,\bar z)\lam(\bar x+\bar y+\bar z,\bar x)\bar J(x,y,z)x\otimes t^{2\bar x+\bar y+\bar z}
$$
$$B'=\lam(\bx,\by)\lam(\bx+\by,\bz)\lam(\bx+\by+\bz,\bx)(\bar J(x,y,z),x)\d_{2\bx+\by+\bz,0}(\bx+\by+\bz),$$
and
$$C'=b'\otimes t^{\bar x+\bar y+\bar z},$$
where
\begin{eqnarray*}
b'&=&-\a' d_x(\bar x+\bar y+\bar z)(xy)z\lam(\bar x,\bar y)\lam(\bar x+\bar y,\bar z)\\
&&+\a' d_x(\bar y)(yz)x\lam(\bar y,\bar z)\lam(\bar y+\bar z,\bar x)\\
&&-\b'  d_x(\bar x+\bar y+\bar z)(zx)y\lam(\bar x,\bar z)\lam(\bar x+\bar z,\bar y)\\
&&-\b'  d_x(\bar z)(zy)x\lam(\bar y,\bar z)\lam(\bar y+\bar z,\bar x)\\
&&-\gamma' d_x(\bar x+\bar y+\bar z)(yz)x\lam(\bar y,\bar z)\lam(\bar y+\bar z,\bar x)\\
&&- \gamma' d_x(\bar y+\bar z)(yz)x\lam(\bar y,\bar z)\lam(\bar y+\bar z,\bar x)\\
&=&-d_x(\bar x+\bar y+\bar z)\big((xy)z+(-1)^{|x||z|}(zx)y+(yz)x\big)\lam(\bar x,\bar y)\lam(\bar x+\bar y,\bar z),\\
\end{eqnarray*}
and
$$\a'=(-1)^{|x||z|},\quad\b'=(-1)^{|y||z|}\andd\gamma'=(-1)^{|x||y|}.
$$

We note that if $d_x\not=0$, then $|x|=0$, in which case
$$y(xz)=(-1)^{|z||x|}(-1)^{|zx||y|}(zx)y=(-1)^{|z||y|}(zx)y,\quad
(xz)y=-(-1)^{|z||x|}(zx)y=-(zx)y.$$
Therefore
\begin{eqnarray*}
b&=&d_x(\bar z)(xy)z\lam(\bar x,\bar y)\lam(\bar x+\bar y,\bar z)\\
&&+\b d_x(\bar y)(zx)y\lam(\bar x,\bar z)\lam(\bar x+\bar z,\bar y)\\
&&+\b  d_x(\bar x+\bar z)(zx)y\lam(\bar x,\bar z)\lam(\bar x+\bar z,\bar y)\\
&&+\b  d_x(\bar z)(zx)y\lam(\bar x,\bar z)\lam(\bar x+\bar z,\bar y)\\
&&+d_x(\bar z)(yz)x\lam(\bar y,\bar z)\lam(\bar y+\bar z,\bar x)\\
&&- \b d_x(\bar x+\bar y+\bar z)(zx)y\lam(\bar x,\bar z)\lam(\bar x+\bar z,\bar y)\\
&=&d_x(\bar z)\big((xy)z+(-1)^{|x||z|}(zx)y+(yz)x\big)\lam(\bar x,\bar y)\lam(\bar x+\bar y,\bar z)\\
&=&d_x(\bar z)\bar J(x,y,z) \lam(\bar x,\bar y)\lam(\bar x+\bar y,\bar z).
\end{eqnarray*}
Similarly, we obtain
\begin{eqnarray*}
b'
&=&-d_x(\bar x+\bar y+\bar z)\big((xy)z+(-1)^{|x||z|}(zx)y+(yz)x\big)\lam(\bar x,\bar y)\lam(\bar x+\bar y,\bar z)\\
&=&-d_x(\bar x+\bar y+\bar z)\bar J(x,y,z) \lam(\bar x,\bar y)\lam(\bar x+\bar y,\bar z) .
\end{eqnarray*}

Now, we investigate under which circumstances,
$$J(\hat x,\hat y,\hat x\hat z)=J(\hat x,\hat y,\hat z)\hat x,$$
for all $\hat x$, $\hat y$, $\hat z$.
Since $\lam$ is a $2$-cocycle
$$\lam(\bar x,\bar y)\lam(\bar x,\bar z)\lam(\bar x+\bar y,\bar x+\bar z)=\lam(\bar x,\bar y)\lam(\bar x+\bar y,\bar z)\lam(\bar x+\bar y+\bar z,\bar x)$$
and
$$\lam(\bar x,\bar y)\lam (\bar x+\bar y,\bar z)=\lam(\bar x,\bar z)\lam (\bar x+\bar z,\bar y)=\lam(\bar y,\bar z)\lam(\bar y+\bar z,\bar x),$$
for all $\bar x,\bar y,\bar z\in G$.
Therefore
$\bar J(\hx,\hy,\hx\hz)=\bar J(\hx,\hy,\hz)\hx$ in $\ll(\aa)$ for all $\hx,\hy,\hz$. Also
$\bar J(\hx,\hy,\hx\hz)=\bar J(\hx,\hy,\hz)\hx$ in $\tilde\ll(\aa)$ if and only if
\begin{equation}\label{nn2}
(\bar J(x,y,z),x)=0\hbox{ for all }x,y,z\in\aa\hbox{ with }|y|+|z|=0.
\end{equation}
In particular, (\ref{nn2}) holds in any anticommutative algebra.
Also $\bar J(\hx,\hy,\hx\hz)=\bar J(\hx,\hy,\hz)\hx$ in $\hat\ll(\aa)$ for all $\hx,\hy,\hz$, if and only if
\begin{equation}\label{nn1}
\begin{array}{c}
\bar J(x,y,x)=0\hbox{ for all }x,y\in\aa,\\
\bar J(x,y,z)=0\hbox{ for all }x,y,z\in\aa\hbox{ with }|x|=0,\\
(\bar J(x,y,z),x)=0\hbox{ for all }x,y,z\in\aa.
\end{array}
\end{equation}
Since $\bar J(x,y,z)$ is invariant under cyclic permutation, the second equality in (\ref{nn1}) is equivalent to $\bar J(x,y,z)=0$ if one of $x$, $y$ or $z$ has degree zero. Now assume $|x|=|y|=|z|=1$. From the first equality in (\ref{nn1}), we have
$\bar J(x+z,y,x+z)=0$. This gives $0=\bar J(x,y,z)+\bar J(z,y,x)=2\bar J(x,y,z)$.
Thus the conditions in (\ref{nn1}) is equivalent to $\aa$ being a Lie (super)algebra.

Finally if in the above computations, we treat as zero all terms containing an element $d\in\v^*$, if follows easily that $\tilde{\ll}(\aa)$ is a Malcev (super)algebra if and only if $(\bar{J}(x,y,z),x)=0$ for all $x,y,z$. In particular, if $\aa$ is a Malcev (super)algebra then so is $\tilde{\ll}(\aa)$. It also follows from our computations that $\aa$ is a Malcev (super)algebra if and only if $\ll(\aa)$ is a Malcev (super)algebra.\qed


\begin{thebibliography}{0}

\bibitem[AABGP]{aabgp}
B. Allison,  S. Azam, S. Berman, Y. Gao, and A.
Pianzola,
  {\it Extended affine {L}ie algebras and their root systems,} {Mem. Amer.
  Math. Soc.} \textbf{126} (1997), no.~603, x+122. 



\bibitem[ABB]{ABB} H. Albuquerque, E. Barreiro, S. Benayadi, {\it  Quadratic Malcev superalgebras with reductive even part,} arXiv:0710.4544v1 [math.RA].

\bibitem[AB]{AB} H. Albuquerquea, S. Benayadi, {\it Quadratic Malcev superalgebras,} J. Pure Appl. Alg. {\bf 187} (3004),
19--45.

%

\bibitem[ABGP]{abgp}
B. Allison,  S. Berman, Y. Gao, and A. Pianzola,
{\it A
  characterization of affine {K}ac-{M}oody {L}ie algebras,} {Comm. Math. Phys.}
  \textbf{185} (1997), no.~3, 671--688. 

\bibitem[AHY]{ahy} S. Azam, R. Hosseini, M. Yousofzadeh, {\it Extended affinization of invariant affine reflection algebras}, Osaka J. Math., {\bf 50} (2013), 1039--1072.

\bibitem[AYY]{ayy} S. Azam, Y. Yoshii, M. Yousofzadeh, {\it Jordan tori for a torsion free abelian group},
\bibitem[CRT]{CRT} H. I. Carrion, M. Rojas, E. Toppan, {\it An $N = 8$ superaffine Malcev algebra and its $N = 8$ Sugawara,} arXiv:hep-th/0105313v1.




\bibitem[E]{E} A. Elduque, {\it On semisimple Malcev algebras,}
Proceed. of AMS, {\bf 107} (1989), 73--82.

\bibitem[GM]{gm} M. Gunaydin, D. Minic, {\it
 Nonassociativity, Malcev Algebras and String Theory 1},
 arXiv:1304.0410 [hep-th].


\bibitem[H-KT]{hkt} R. H\o egh-Krohn, B. Torr\'esani, {\it Classification and construction of quasi-simple Lie algebras},  J. Funct. Anal. 89 (1990), 106?136.

\bibitem[K]{kac}
V. Kac,  {Infinite-dimensional {L}ie algebras,} third
ed., Cambridge
  University Press, Cambridge, 1990. 


\bibitem[LN]{nllfrs}
O. Loos,   and E. Neher, {\it Locally finite root systems,}
{Mem. Amer.
  Math. Soc.} \textbf{171} (2004), no.~811, x+214. 


\bibitem[MY]{my}
J. Morita,  and Y. Yoshii, {\it Locally extended affine {L}ie
algebras,} {J.
  Algebra} \textbf{301} (2006), no.~1, 59--81. 

  \bibitem[Nee]{nee} K-H. Neeb, {\it Integrable roots in split graded Lie algebras,} J. Algebra, {\bf 225} (2000), 534--580.

  \bibitem[NS]{NS} K-H. Neeb, N. Stumme, {\it The classification of locally finite split simple Lie algebars}, J. reine angew. Math, (2001), 25-53.

\bibitem[Neh]{neh}
E. Neher {\it Extended affine {L}ie algebras and other
 generalizations of affine {L}ie algebras - a survey,} ArXiv e-prints (2008).


\bibitem[O]{O} E. P. Osipov, {\it
Kac-Moody-Malcev and super-Kac-Moody-Malcev algebras,} Phisics Letter B, {\bf 274} (1992), 341--344.

\bibitem[Sa1]{sa1} A. A. Sagle, {\it Malcev algebras},

\bibitem[Sa2]{sa2} A. Sagle, {\it On simple extended Lie algebras over fields of  characteristic zero,} Pacific J. Math,
{\bf 15} (1965), 621-648.

\bibitem[Stu]{stu} N. Stumme, {\it The Structure of locally finite split Lie algebras,} J.  Algebra {\bf 220}, (1999) 664?693.

\bibitem[Yos1]{yos1} Y. Yoshii, {\it Lie tori- A simple characterization of extended affine Lie algebras} Publ. RIMS, Kyoto Univ.
42 (2006), 739�762.

\bibitem[Yos2]{yos2} Y. Yoshii, {\it Lie $G$-tori of symplectic type} Quart. J. Math. 57 (2006), 425�448.


\bibitem[You1]{you1}
M. Yousofzadeh, {\it Extended affine Lie superalgebras,} arXiv:1309.3766v3 [math.QA].

\bibitem[You2]{you2}
M. Yousofzadeh, {\it Extended affine root supersystems,} arXiv:



\end{thebibliography}
\end{document}